\documentclass{svproc}

\usepackage{amsmath,amssymb,graphics,epsfig,color,enumerate,mathrsfs}

\usepackage[showlabels,sections,floats,textmath,displaymath]{preview}

\newtheorem{Theorem}{Theorem}

\newtheorem{Prop}{Proposition}

\usepackage{dsfont} 
\definecolor{refkey}{gray}{.75}
\definecolor{labelkey}{gray}{.5}

\newcommand{\bbE}{{\ensuremath{\mathbb E}} }

\newcommand{\bbL}{{\ensuremath{\mathbb L}} }

\newcommand{\bbP}{{\ensuremath{\mathbb P}} }
\newcommand{\bbQ}{{\ensuremath{\mathbb Q}} }
\newcommand{\bbR}{{\ensuremath{\mathbb R}} }

\newcommand{\bbZ}{{\ensuremath{\mathbb Z}} }

 \let\b=\beta     \let\e=\varepsilon
        \let\l=\lambda
      \let\o=\omega      
   \let\t=\tau

\let\O=\Omega      
\usepackage{url}

\begin{document}
\mainmatter              
\title{1D Mott variable-range hopping with external field}
\titlerunning{1D Mott variable-range hopping}  
%
\author{Alessandra Faggionato}
%
\authorrunning{Alessandra Faggionato} 
%
\institute{University La Sapienza, P.le Aldo Moro 2, Rome, Italy,\\
\email{faggiona@mat.uniroma1.it}}
\maketitle

\begin{abstract}
 Mott variable-range hopping is  a fundamental mechanism for electron transport in disordered solids in the regime of strong Anderson localization. We give a brief description of 
  this mechanism, 
 recall  some results concerning the behavior of the conductivity at low temperature and   describe in  more   detail 
recent results (obtained in collaboration with N. Gantert and M. Salvi) concerning the one-dimensional Mott variable-range hopping 
under  an external field\footnote{If you want to cite this proceeding and in particular  some result contained in it, please cite also
the article where the result appeared, so that the contribution of my coworkers can be recognised.  This can be important  for bibliometric reasons.}
.
\end{abstract}

\section{ Mott variable-range hopping }

Mott  variable range hopping is a mechanism of phonon-assisted electron transport taking place in amorphous solids (as doped semiconductors)  in the regime of strong Anderson localization. It has been introduced by N.F.  Mott in order to explain the anomalous 
non-Arrhenius decay of the conductivity at low temperature \cite{M1,M3,MD,POF,SE}.

Let us consider a  doped semiconductor, which is given by a semiconductor
with randomly located foreign atoms (called \emph{impurities}).  We write $\xi:=\{x_i\}$ for the set of impurity sites. For simplicity we treat spinless electrons. Then, 
due to  Anderson localization, a generic  conduction electron  is  described by a quantum  vawefunction  localized around some impurity site $x_i$, whose energy is denoted by $E_i$.  This allows, at a first approximation, to think of the conduction  electrons as classical particles which can lie only on the impurity sites, subject  to the constraint of site exclusion (due to Pauli's exclusion principle). As a consequence, a microscopic configuration is described by an element $\eta\in \{0,1\}^{\xi}$,  where $\eta_{x_i}=1$ if and only if an electron is localized around the impurity site $x_i$. The dynamics  is then described by an exclusion process, where  the 
 probability  rate  for a jump from $x_i$ to $x_j$ 
  is given by (cf. \cite{AHL})
\begin{equation}
 \mathds{1}(\eta_{x_i}=1,\; \eta_{x_j}=0) \exp\{-2\zeta  |x_i-x_j| - \beta \{ E_j -E_i \}_+\} 
 \,.
 \label{nove}
\end{equation} 
Above  $\b=1/kT$ ($k$ being the Boltzmann's constant and $T$ the absolute temperature), while $1/\zeta$ is the localization length. A physical analysis suggests that the energy marks $\{E_i\}$ can be modeled 
by i.i.d. random variables. In inorganic doped semiconductors, when the Fermi energy is set equal to zero,  the common distribution $\nu$ of $E_i$ is of the form 
$c\, |E|^{\alpha} dE$ on some interval $[-A,A]$, where $c$ is  the normalization constant and $\alpha$ is a nonnegative exponent.

At  low temperature (i.e. large $\beta$)  the form of the jump rates (\ref{nove}) suggests that long jumps can be facilitated when the  energetic cost  $ \{ E_j -E_i \}_+$ is small. This facilitation  leads to an anomalous conductivity behavior for $d\geq 2$. Indeed, according to Mott's law, in an isotropic medium the conductivity matrix $\sigma (\beta)$ can be approximated (at logarithmic scale) by  $ \exp 
 \bigl( - c\, \beta^{\frac{\alpha+1}{\alpha+1+d} }\bigr)\mathds{1}$, where $\mathds{1}$ denotes the identity matrix  and $c$ is a suitable  positive constant $c$ with negligible temperature-dependence.

The mathematical analysis of the above    exclusion 
process presents several technical challenges and has been performed only when  $\xi\equiv \bbZ^d$ (absence of geometric disorder) and with jumps restricted to nearest-neighbors (cf. \cite{FMa,Q}). This last assumption does not fit with the low temperature regime, where  anomalous conductivity takes place. 

To investigate Mott  variable range hopping at low temperature, 
in the regime of low impurity density some effective models  have been proposed. One is given by the 
random Miller-Abrahams resistor network \cite{AHL,MA}. Another effective model is the following  (cf. \cite{FSS}): one 
 approximates  the localized electrons  by classical non--interacting (independent) particles moving according to random walks with  jump  probability  rate 
given by the  transition rate (\ref{nove}) multiplied  by a suitable factor which keeps  trace of the exclusion principle. To be more precise, 
given $\gamma \in \bbR$, we call  $\mu_\gamma$  the product probability measure on $\{0,1\}^\xi$ with 
$\mu_\gamma (\eta_{x_i})= \frac{e^{-\beta (E_i-\gamma)} }{ 1+ e^{-\beta (E_i-\gamma)}  }
$. Then it is simple to check that $\mu_\gamma$ is a reversible distribution for the  exclusion process.
In the independent particles approximation, the probability rate for a jump from $x_i$ to $x_j$ is given by 
\begin{equation}\label{fisarmonica}
\mu_\gamma (\eta_{x_i}=1, \; \eta_{x_j} =0)  \exp\{-2 \zeta  |x_i-x_j| - \beta \{ E_j -E_i \}_+\}\,.
\end{equation} 
It is simple to check that at low temperature, i.e. large $\beta$, (\ref{fisarmonica}) is well approximated by (cf. \cite{AHL})
\begin{equation}\label{viaggio}
 \exp\{  - 2 \zeta |x_i-x_j| - \frac{\beta}{2} ( |E_{i} -\gamma |+|E_{j} -\gamma | + |E_{i}-E_{j}|)\}\,. 
 \end{equation}
 In what follows, without loss of generality,  we shift the  energy so that the Fermi energy $\gamma$ equals zero, we take $2\zeta=1$ and we replace $\beta/2$ by $\beta$. Since the above random walks are independent, we can restrict to the analysis of a single random walk, which we call  \emph{Mott  random walk}.  
 
%
 
 \section{Mott  random walk and bounds on the diffusion matrix}
We give the formal definition of  Mott  random   
walk  as  random walk in a random environment. The environment is given by 
a marked simple point process $\o= \{ (x_i,E_i)\}
$ where $\{x_i\}\subset \bbR^d$ and the (energy) marks are i.i.d. random variables with common distribution $\nu$ having support on some finite interval $[-A,A]$.
 Given a realization of the environment $\o$,  Mott random walk is  the continuous-time random walk $X_t^{\o}$
with state space    $ \{x_i\}$ and probability rate for a jump from $x_i$ to $x_j \not=x_i$ given by 
\begin{equation}\label{germania}  r_{x_i,
x_j}(\o) = \exp \left\{
 -|x_i-x_j| -\beta (|E_i|+|E_j|+|E_i-E_j|)\right\} \,.
\end{equation}
As already mentioned, one expects for $d\geq 2$ that the contribution to the transport of long jumps dominates as $\beta \to \infty$.  In \cite{CFP}  a quenched invariance principle for Mott  random walk has been proved. Calling $D(\beta)$ the diffusion matrix of the limiting Brownian motion, in  \cite{FM,FSS} bounds in agreement with Mott's law have been obtained for the diffusion matrix. More precisely, under very general conditions on the isotropic environment  and taking $\nu$ of the form $c\, |E|^\alpha dE$ for some $\alpha \geq 0$ and on some interval $[-A,A]$,  
it has been proved that for suitable $\beta$-independent  positive constant $c_1, c_2,\kappa_1, \kappa_2$ it holds 
\begin{equation}\label{sigari}
c_1 \exp \left\{ - \kappa_1\,\beta^{\frac{\alpha+1}{\alpha+1+d}}   \right\} \mathds{1}  \leq D(\beta)\leq c_2 \exp \left\{
- \kappa_2\, \beta^{\frac{\alpha+1}{\alpha+1+d}  }
\right\} \mathds{1} \,.
\end{equation}
We point out that for a genuinely nearest--neighbor random walk $D(\beta)$ would have an Arrenhius decay, i.e.  $D(\beta) \approx e^{-c \beta} \mathds{1}$, thus implying that  (\ref{sigari}) is determined by  long jumps.

In dimension $d=1$ long jumps do not dominate.  Indeed,  the following has been derived for $d=1$ in \cite{CF} when $\o$ is a  renewal marked simple point process. We label the points in increasing order, i.e. $x_i< x_{i+1}$,  take  $x_0=0$ and assume that $\bbE[ x_1^2]<\infty$. Then the following holds \cite{CF}:

\begin{itemize}
\item[(i)]  If   $\bbE[ e^{x_1}] <\infty$, then a quenched invariance principle holds and the diffusion coefficient satisfies 
\begin{equation}\label{sigari2}
c_1 \exp \left\{ - \kappa_1\,\beta   \right\}    \leq D(\beta)\leq c_2 \exp \left\{
- \kappa_2\, \beta   
\right\} \,,
\end{equation}
for $\beta$--independent positive constants $c_1, c_2,\kappa_1, \kappa_2$;
\item[(ii)]  If    $\bbE[ e^{x_1}] =\infty$, then  the random walk is subdiffusive. More precisely, an  annealed invariance principle holds with zero diffusion coefficient.\end{itemize}

We point out that the above 1d results hold also for a larger class of   jump rates and energy distributions $\nu$ \cite{CF}.
Moreover, we stress  that the above bounds (\ref{sigari}) and (\ref{sigari2}) refer to the diffusion matrix $D(\beta)$ and not to the conductivity matrix $\sigma(\beta)$. On the other hand, believing in the  Einstein relation (which states that $\sigma(\beta)=\beta D(\beta)$), the above bounds would extend to $\sigma (\beta)$, hence 
one would recover lower and upper bounds on the conductivity matrix in agreement with  the physical Mott  law for $d\geq 2$ and with the Arrenhius--type decay for $d=1$. 
The rigorous derivation of the Einstein relation for Markov processes in random environment is in general a difficult task and has been the object of much investigation also in the last years (cf. e.g. \cite{GMP,GGN,KO1,KO2,LR,L1,L2,MP2}). 
In what follows, we will concentrate on the effect of perturbing 1d  Mott  random walk by an external field and on the validity of the Einstein relation. 
 Hopefully, progresses on the Einstein relation for Mott random walk in higher dimension will be obtained in the future. 

 \section{Biased 1d Mott random walk}
We take $d=1$ and  label points $\{x_i\}$ in increasing order with the convention that $x_0:=0$ (in particular, we assume the origin to be an impurity site).
  It is convenient to define $Z_i$ as the interpoint distance $Z_i:= x_{i+1}-x_i$.

%

\begin{figure}
 \includegraphics[width=10cm]{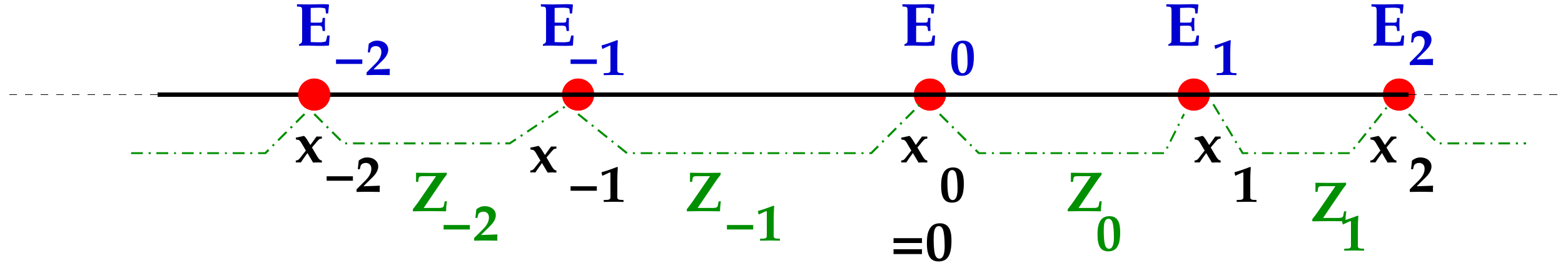}
\caption{Points $x_i$, energy marks $E_i$ and interpoint distances $Z_i$. }
\end{figure}

We make the following assumptions:

\begin{itemize}
\item[(A1)] The  random sequence $( Z_k,E_k)_{k \in \bbZ}$ is  stationary and ergodic  w.r.t. shifts;
\item[(A2)]   $\bbE[Z_0] $ is finite;
\item[(A3)]   $\bbP ( \o = \t_\ell \o)$ is zero for  all  $\ell \in \bbZ\setminus \{0\}$;   
\item[(A4)]    There exists  some constant $  d>0$ satisfying
 $\bbbp(Z_0 \geq d)=1$.
\end{itemize}

Note that we do not restrict to the physically relevant energy mark distributions $\nu$ which are  of the form $c\,|E|^\alpha dE$ on  some interval $[-A,A]$.

Given $\l\in [0,1)$ we consider the biased 
generalized Mott  random walk  $X_t^{\o,\l}$ on $\{x_i\}$ with jump probability rates given by 
\begin{equation}\label{def_rates}
r_{x_i,x_j} ^\l(\o)= \exp \left\{- |x_i-x_j|  +\l(x_j-x_i)   - u(E_i,E_j)  \right\}\,, \qquad x_i \not =x_j\,,
\end{equation}
and starting at the origin. 
Above,  $u$ is a given bounded and symmetric function. Note that we do not restrict to (\ref{germania}). The special form \eqref{germania} is relevant when studying the regime $\b \to \infty$, on the other hand here we are interested in the system at a fixed temperature (which is included in the function $u$) under the effect of an external field. 

 In the rest, it is convenient to set $r_{x,x}^\l(\o)\equiv 0$. We also point out that one can easily prove that the random walk 
 $X_t^{\o,\l}$  is well defined since $\l\in[0,1)$.
 \smallskip
 
The following result, obtained in \cite{FGS1}, concerns the ballistic/sub--ballistic regime:
 
\begin{Prop}{\sc \cite{FGS1}} For $\bbP$--a.a. $\o$ the random walk  $X_t^{\o ,\l}$ is transient to the right, i.e. $\lim _{t \to \infty}   X_t^{\o ,\l}=+\infty$  a.s. 
\end{Prop}

\begin{Theorem}\label{30minuti} {\sc \cite{FGS1}}
Fix $\l\in (0,1)$.
\begin{itemize}

\item[(i)] If  $\bbE \big[ e ^{ (1-\l) Z_0}\bigr] < \infty$ and $u$ is continuous, then for $\bbP$--a.a. $\o$ the following limit exists 
\[v_X(\l):=  \lim _{t \to \infty} \frac{ X_t^{\o ,\l} }{ t}  \qquad a.s.\]
and moreover it is deterministic, finite and strictly positive.

\item[(ii)]   If  $\bbE \big[ e^{ -(1+\l)Z_{-1} + (1-\l) Z_0}\bigr] =\infty$, then for $\bbP$--a.a. $\o$ it holds
\[  v_X(\l):= \lim _{t \to \infty} \frac{ X_t^{\xi ,\l} }{ t} =0\qquad a.s.\]
\end{itemize}
\end{Theorem}

As discussed in \cite{FGS1}, the condition $\bbE \big[ e ^{ (1-\l) Z_0}\bigr] =\infty$ does not imply that $v_X(\l)=0$.
On the other hand,
 if 
  $(Z_k)_{ k \in \bbZ}$ are i.i.d. (or even if $Z_k,Z_{k+1}$ are independent for every $k$) and $u$ is continuous,
  then the above two cases (i) and (ii) in Theorem \ref{30minuti} are exhaustive and one concludes that 
  $ \bbE \big[ e ^{ (1-\l) Z_0}\bigr]<\infty$ if and only if  $ v_X(\l) > 0$, otherwise $v_X(\l)=0$.  
  
  \medskip

  The above Theorem \ref{30minuti}   extends also to the jump process associated to      $X^{\o,\l}_t$, i.e. to the discrete time random walk 
 $Y^{\o,\l} _n$ with probability $p_{x_i,x_k}^\l(\o)$  of a jump from $x_i$ to $x_k \not = x_i$ given by 
 \[
 p_{x_i,x_k}^\l(\o)= \frac{ r_{x_i,x_j} ^\l(\o)  }{  \sum_k r_{x_i,x_k} ^\l(\o)   }\,.
 \]
  In particular,   if $\bbE \big[ e ^{ (1-\l) Z_0}\bigr] < \infty$ and $u$ is continuous then the random walk $Y^{\o,\l} _n$ is ballistic  ($v_Y(\l)>0$), while if  $\bbE \big[ e^{ -(1+\l)Z_{-1} + (1-\l) Z_0}\bigr] =\infty$ then 
 the random walk $Y^{\o,\l} _n$ 
    is sub--ballistic (i.e. $v_Y(\l)=0$).
   We point out  that  indeed the result has been proved in \cite{FGS1} first for the random walk $Y^{\o,\l} _n$ and then extended to the continuous time case by a random time change argument.

\smallskip

We write $\t_{x} \o$ for the environment $\o$ translated by $x\in \bbR$, more precisely we set 
$\t_{x} \o:= \{ (x_j-x, E_j) \} $  if $\o=\{(x_j,E_j)\}$. We recall that the environment viewed from the walker  $Y^{\o,\l} _n$  is given by the discrete time Markov chain $( \t_{ Y^{\o,\l} _n} \o)_{n\geq 0}$. This is the crucial object to analyze in the ballistic regime. 
The following result concerning  the environment viewed from the walker  $Y^{\o,\l} _n$
is indeed at the basis of the  derivation of Theorem \ref{30minuti}--(i) as well as the starting point for the analysis of the Einstein relation.
 \begin{Theorem}\label{1ora} {\sc \cite{FGS1,FGS2}} Fix $\l\in(0,1)$.  
 Suppose   that   $\bbE \big[ e ^{ (1-\l) Z_0}\bigr] < \infty$  and that  $u$ is  continuous.
Then the  environment viewed from  the walker $Y^{\o,\l} _n$ admits  an invariant and ergodic distribution  $\bbQ_\l$ mutually absolutely continuous w.r.t.  $\bbP$. Moreover, it holds
\[   v_Y(\l)=  \bbQ_\l \bigl[  \varphi_\l \bigr] \;\; \text{ and }\;\;  v_X(\l)=\frac{  v_Y(\l)  }{  \bbQ_\l\Big[ 1/ ( \sum _k r_{0,x_k} ^\l  ) \Big]}\,, \]
where $\varphi_\l$ denotes the local drift, i.e.
$\varphi_\l (\o):= \sum _{i} x_i p_{0, x_i}^\l(\o) $.
\end{Theorem}

We point that  for $\l=0$ the probability distribution $\bbQ_0$ defined as
\[ d \bbQ_0= \frac{  \sum _k r^{\l=0}_{0,x_k}  }{ \bbE [   \sum _k r^{\l=0}_{0,x_k} ]}d\bbP\]
is indeed reversible for the environment viewed from the walker $Y^{\o,\l=0} _n$, and that 
 $v_Y(0)=v_X(0)=0$. In what follows, when $\l=0$ we will often drop $\l$ from the notation (in particular, we will write simply $ r_{x_i,x_j} (\o) $,  $p_{x_i,x_k}(\o)$, $X^{\o }_t$,  $Y^{\o }_n$).

\bigskip

The proof of Theorem \ref{1ora} takes inspiration from the paper \cite{CP}. The main technical difficulty comes from the presence of arbitrarily long jumps, which does not allow to use  standard  techniques based on regeneration times. Following the method developed in \cite{CP} we have considered, for each positive integer $\rho$, the random walk   
 obtained from  $Y^{\o ,\l}_n$ by suppressing  jumps   between sites $x_i,x_j$ with $|i-j|>\rho$ \cite{FGS1}. For this $\rho$-indexed random walk,  under the same hypothesis of Theorem \ref{1ora},    we have proved  that there exists a distribution 
$\bbQ_\l ^{(\rho)}$ which is   invariant and ergodic  for the associated environment viewed from the walker, and that  $\bbQ_\l ^{(\rho)}$ is  mutually absolutely continuous w.r.t.  $\bbP$. In particular, the methods developed in \cite{CP}  provide a probabilistic representation of the Radon--Nykodim $\frac{d\bbQ_\l^{(\rho)}}{d \bbP}$,  which (together with a suitable analysis  based on potential theory) allows to prove that 
 $\bbQ_\l ^{(\rho)}$ weakly converges to  $\bbQ_\l$, and that $\bbQ_\l$ has indeed the nice properties stated in Theorem \ref{1ora}.
%
%
 \section{Linear response and Einstein relation for the biased 1d Mott random walk}
 
 Let us assume again (A1), (A2), (A3), (A4) as in the previous section.
 The probabilistic representation of  the Radon--Nykodim derivative  $\frac{d\bbQ_\l^{(\rho)}}{d \bbP}$ mentioned in the above section is the starting point for 
the derivation of  estimates on the Radon--Nykodim derivative  $\frac{d \bbQ_\l}{d \bbQ_0 }$:
\begin{Prop}\label{andrea} {\sc \cite{FGS2}}  Suppose that for some   $p\geq 2$  it holds $\bbE\bigl[ e ^{p Z_0 }\bigr] < +\infty $. Fix $\l_0 \in (0,1)$.  Then  
 \[
\sup _{\l \in [0,\l_0]} \Big \| \frac{d \bbQ_\l}{d \bbQ_0 }\Big \| _{L^p (\bbQ_0) }  <\infty\,.
 \]
\end{Prop}
The above proposition  allows to prove the 
 continuity of the expected value  $\bbQ_\l(f) $ of suitable functions $f$:
\begin{Theorem}\label{teo_continuo}  {\sc \cite{FGS2}} 
Suppose that $\bbE[{\rm e}^{p Z_0}] <\infty$  for some $p\geq 2$ and let $q$ be the conjugate exponent of $p$, i.e.~$q$ satisfies  $\frac{1}{p}+\frac{1}{q}=1$. 
Then,  for any $f \in L^q(\bbQ_0)$ and $\l \in [0,1)$,   it holds that  $f \in L^1(\bbQ_\l)$ and    the map 
 \begin{equation}\label{conti}
[0,1) \ni \l \mapsto  \bbQ_\l(f) \in \bbR 
\end{equation}
is continuous.
\end{Theorem}
Without entering into the details  of the proof (which can be found in \cite{FGS2}) we give some comments on the derivation of Theorem 
\ref{teo_continuo} from Proposition \ref{andrea}. To this aim, we take for simplicity $p=2$.  Hence we are supposing that 
 $\bbE({\rm e}^{2 Z_0}) <\infty$ and that $f \in L^2(\bbQ_0)$, and we want to prove that $  f \in L^1(\bbQ_\l)$ and that the function in \eqref{conti} is continuous. For simplicity, let us restrict to its continuity at $\l=0$.
 The fact that $  f \in L^1(\bbQ_\l)$ follows by writing $\bbQ_\l(f) = \bbQ_0(  \frac{d \bbQ_\l}{d\bbQ_0} f )$ and then by 
 applying Schwarz inequality  and Proposition \ref{andrea}.  The proof of the continuity at $\l=0$ is more involved.
We recall that by Kakutani's  theorem balls are compact for the $L^2(\bbQ_0)$--weak topology. Hence, due to Proposition \ref{andrea}, the family  of Radon--Nykodim derivatives $\frac{d \bbQ_\l}{d \bbQ_0 }$, $\l \in [0,\l_0]$, is  relatively compact for the  $L^2(\bbQ_0)$--weak topology. In \cite{FGS2} we then prove that any limit point of this family is given by $\mathds{1}$. As a byproduct of the representation    $\bbQ_\l(f) = \bbQ_0(  \frac{d \bbQ_\l}{d\bbQ_0} f )$ and of the weak convergence  $\frac{d \bbQ_\l}{d \bbQ_0 } \rightharpoonup  \mathds{1}$, we get the continuity of \eqref{conti} at $\l=0$.

\bigskip

We now move to the study of 
$\partial_{\l=0} \bbQ_\l(f)$. To this aim we introduce the operator    $\bbL_0 : L^2(\bbQ_0) \to L^2(\bbQ_0)$ as 
\[ \bbL_0 f(\o)=\sum _k p_{0,x_k} (\o) [ f(\t_{x_k} \o )- f(\o) ]\,, \qquad f\in L^2(\bbQ_0)\,.
\]
We recall that a function $f$  belongs to 
$L^2(\bbQ_0)\cap H_{-1}$ if there exists $C>0$ such that 
\[ 
|\langle f, g \rangle  | \leq C \langle g, -\bbL_0  g \rangle^{1/2} \qquad  \forall g\in L^2(\bbQ_0)\,. 
\]
Above $\langle \cdot, \cdot \rangle$ is the scalar product in $L^2(\bbQ_0)$.
Due to the theory developed by Kipnis and Varadhan \cite{KV}, for any $f\in  L^2(\bbQ_0)\cap H_{-1}$, 
 we have the weak convergence 
\[ \frac{1}{\sqrt{n} } \bigl( \sum_{j=0}^{n-1}  f(\o_j), \sum_{j=0}^{n-1} \varphi (\o_j) \bigr) \stackrel{n\to \infty }{\rightarrow } (N^f , N^\varphi)\]
for a suitable  2d gaussian vector $(N^f , N^\varphi)$. Above $\varphi$ denotes the local drift $\varphi_\l$ with $\l=0$ (cf. Theorem \ref{1ora}) and $\o_j=\t_{Y_j^\o} \o$, i.e. $(\o_n)_n$ represents the environment viewed from the walker $Y_n^\o$.

Finally, we need another ingredient coming from the theory of square integrable forms  in order to present our next theorem. 
 We consider the space $\O\times \bbZ$ endowed with the  
measure $M$ defined by 
 \[ 
 M(v) = \bbQ_0 \Big[ \sum _{k\in\bbZ} p_{0,{x_k}} v(\cdot, k ) \Big] \,, \qquad \forall v :\O \times \bbZ\to \bbR \qquad \text{Borel, bounded}
 \,.\]
 A generic Borel function $v: \O \times \bbZ \to \bbR$ will be called a \emph{form}.   $L^2 ( M)$   is known as the space of \emph{square integrable forms}.  Given a  function $g =g(\o) $ we define 
 \begin{equation}\label{cirm}
  \nabla g (\o, k):= g(\t_{x_k} \o) - g(\o) \,.
 \end{equation}
If  $g \in L^2(\bbQ_0) $, then  $\nabla g\in L^2(M)$. The closure in $M$  of the subspace $\{ \nabla g \,:\, g\in L^2(\bbQ_0) \}$ is the set of   the so called \emph{potential forms} (its orthogonal subspace is given by the so called \emph{solenoidal forms}).

  Take again $f \in H_{-1}\cap L^2(\bbQ_0)$ and, given $\e>0$, define $g^{f} _\e \in L^2(\bbQ_0)$ as the unique solution of the equation
 \begin{equation}\label{def_gigi}
 (\e -\bbL_0) g_{\e }^{f}= f \,.
 \end{equation}
 As discussed with  more details in \cite{FGS2}   as $\e$ goes to zero the family of potential forms  $\nabla g_\e^{f}$ converges in $L^2(M)$ to a potential form $h^{f}$:
 \begin{equation}\label{hf} 
 h^{f} = \lim _{\e \downarrow 0 } \nabla g_\e ^{f} \qquad \text{ in } L^2(M)\,.
 \end{equation}

\begin{Theorem}\label{teo_derivo}$ $
 {\sc \cite{FGS2}}
Suppose $\bbE({\rm e}^{p Z_0}) <\infty$  for some $p>2$.  Then, for any $f \in H_{-1}\cap L^2(\bbQ_0)$,  $  \partial_{\l=0} \bbQ_\l (f) $  exists.
Moreover the following two probabilistic representations hold with $h=h^f$ (see \eqref{hf}):
 \begin{equation}\label{schorle}
 \partial_{\l=0} \bbQ_\l (f)   = \begin{cases} \bbQ_0 \bigl[ \sum_{k\in \bbZ} p_{0,x_k} (x_k -\varphi) h(\cdot, k) \bigr]\\-{\rm Cov}(N^f, N^\varphi)
\end{cases}\,.
\end{equation}
\end{Theorem}
We point out that a covariance representation of $  \partial_{\l=0} \bbQ_\l (f) $ as the second one in \eqref{schorle} appears also
in \cite{GGN} and \cite{MP2}.

\medskip
We give some comments on the derivation of Theorem \ref{teo_derivo} up to the first representation in \eqref{schorle}.
Fix $f \in H_{-1}\cap L^2(\bbQ_0)$, thus implying that $\bbQ_0(f)=0$. Given $\e>0$ take $g_\e\in L^2(\bbQ_0)$ as the unique solution of the equation
 $(\e -\bbL_0) g_{\e }= f$, i.e. $g_\e=g_\e^f$ with $g_\e^f$ as in \eqref{def_gigi}. Then we can write 
\begin{equation}\label{bagno}
\frac{\bbQ_\l (f) -\bbQ_0(f)}{\l} =\frac{\bbQ_\l (f) }{\l} =  \frac{ \e \bbQ_\l (g_\e)}{\l}- \frac{ \bbQ_\l ( \bbL_0 g_\e)}{\l} \,.
\end{equation}
The idea is to 
take first the limit $\e \to 0$, afterwards the limit $\l \to 0$. By  the results  of 
Kipnis and Varadhan \cite{KV} we have that  $\e \bbQ_\l (g_\e)$ is negligible as $\e\to 0$. 
 Hence we have $\partial_{\l=0} \bbQ_\l (f)  = - \lim _{\l \to 0}  \frac{ \bbQ_\l ( \bbL_0 g_\e)}{\l}  $ if the latter exists. On the other hand we have the following identity and approximations for $\e,\l$ small:
\begin{equation}\label{urla}
\begin{split}
- \frac{ \bbQ_\l [\bbL_0 g_\e]}{\l} &=   \bbQ_\l \Big[ \frac{(\bbL_\l- \bbL_0 ) g_\e}{\l}\Big]  =  {\bbQ_\l} \Big[     \sum _{k \in \bbZ } \frac{p^\l_{0,x_k} -  p_{0,x_k} }{\l}( g_\e (\t_{x_k} \cdot )- g_\e) \Big] 
\\&
\approx   {\bbQ_\l} \Big[     \sum _{k \in \bbZ } \partial_{\l=0} p^\l_{0,x_k}  \,  h(\cdot,k) \Big]  \approx   {\bbQ_0} \Big[     \sum _{k \in \bbZ } \partial_{\l=0} p^\l_{0,x_k}\,   h(\cdot, k) \Big]\\
&=   \bbQ_0 \bigl[ \sum_{k\in \bbZ} p_{0,x_k} (x_k -\varphi) h(\cdot, k)  \bigr]\,,
\end{split}
\end{equation}
where $ \bbL_\l f(\o)=\sum _k p^\l_{0,x_k} (\o) [ f(\t_{x_k} \o )- f(\o) ]$ and $h=h^f$ (see \eqref{hf}).
Roughly,  the first identity  in \eqref{urla} follows from the stationary of $\bbQ_\l$ for the environment viewed from $Y_n^{\o,\l}$,  the first approximation in the second line follows from \eqref{hf},  the second  approximation in the second line follows from Theorem \ref{teo_continuo}, the identity in the third line follows from the equality $ \partial_{\l=0} p^\l_{0,k} =p_{0,x_k} (x_k -\varphi)$.

The above steps are indeed rigorously proved in \cite{FGS2}. 
Since, as already observed, $\partial_{\l=0} \bbQ_\l (f)  = - \lim _{\l \to 0}  \frac{ \bbQ_\l ( \bbL_0 g_\e)}{\l}  $, the content of Theorem \ref{teo_derivo} up to the first representation in \eqref{schorle} follows from \eqref{urla}.

 %
%
%
%

\medskip
We conclude this section with the Einstein relation. To this aim we denote by 
 $D_X$ the  diffusion coefficient associated to $X^\o _t$   and by $D_Y$ the  diffusion coefficient associated to  $Y^{\o} _n$. 
$D_X$ is the variance of   the Brownian motion to which  $X^\o_t$ converges under diffusive rescaling, and  a similar
definition holds 
 for $D_Y$.

\begin{Theorem}\label{sciopero} $ $
  {\sc \cite{FGS2}}
  The following holds:
\begin{itemize}
\item[(i)] If  $\bbE[{\rm e}^{2 Z_0}] <\infty$, then  $v_Y(\l)$ and   $v_{X}(\l)$ are continuous functions of $\l$;
\item[(ii)] If  $\bbE[{\rm e}^{p Z_0}] <\infty$  for some $p>2$, then the Einstein relation is fulfilled, i.e.
\begin{equation}\label{einstein}
\partial_{\l=0} v_Y (\l) = D_Y\qquad \text{ and } \qquad \partial _{\l=0}v_{X} (\l) = D_X\,.
\end{equation}
\end{itemize}
\end{Theorem}
If we make explicit the temperature dependence in the jump rates \eqref{def_rates} we would have 
\[
r_{x_i,x_j} ^\l(\o)= \exp \left\{- |x_i-x_j|  +\l\b(x_j-x_i)   - \b u(E_i,E_j)  \right\}\,,
\]
where $\l$ is the strength of the external field. 
Then  Einstein  relation \eqref{einstein}  takes the more familiar (from a physical viewpoint) form
\[
\partial_{\l=0} v_Y (\l, \b)= \b   D_Y(\b) \qquad \text{ and }\qquad  \partial _{\l=0}v_{X} (\l,\b) = \b D_X(\b)\,.
\]

We conclude by giving some ideas behind the proof of the Einstein relation for $v_Y(\l)$ in Theorem \ref{sciopero}. We do not fix the details, but only the main arguments. By applying  Theorem \ref{teo_derivo} with $f:=\varphi$ and using the first representation in the r.h.s. of \eqref{schorle}, one  can express $ \partial_{\l=0} \bbQ_\l[ \varphi] $ as a function of $h=h^\varphi$:
\[
 \partial_{\l=0} \bbQ_\l[ \varphi] = \bbQ_0 \Bigl[ \sum_{k\in\bbZ} p_{0,x_k} (x_k -\varphi) h(\cdot, k) \Bigr]\,.
 \]
 Recall that  $v_Y(\l)=  \bbQ_\l \bigl[  \varphi_\l \bigr] $ (cf. Theorem \ref{1ora})  and that $v_Y(0)=0$. In \cite{FGS2} we have proved the following approximations and identities (for the last identity see Section 9 in \cite{FGS2}) :
\begin{equation*}
\begin{split}
\frac{ v_Y(\l) -v_Y(0)}{\l}&= \frac{ v_Y(\l) }{\l}= \frac{ \bbQ_\l [\varphi_\l] }{\l} =  \bbQ_\l   \bigl[ \frac{ \varphi_\l-\varphi }{\l}\bigr]+\frac{ \bbQ_\l[ \varphi]-\bbQ_0[\varphi] }{\l}\\
&\approx \bbQ_0 \bigl[ \partial_{\l=0} \varphi_\l] + \partial_{\l=0} \bbQ_\l[ \varphi]\\
& = \bbQ_0 \bigl[ \partial_{\l=0} \varphi_\l] +  \bbQ_0 \Bigl[ \sum_{k\in\bbZ} p_{0,x_k} (x_k -\varphi) h(\cdot,k) \Bigr]\\
& =\bbQ_0 \Bigl[ \sum_{k\in\bbZ} p_{0,x_k} (x_k -\varphi)(x_k+ h(\cdot,k)) \Bigr]  = D_Y\,.
\end{split}
\end{equation*}

\bigskip

\bigskip

\noindent
{\bf Acknowledgements}.  It is a pleasure to thank 
 the Institut Henri Poincar\'e  and the Centre Emile Borel  for the kind hospitality and support during the trimester  ``Stochastic Dynamics Out of Equilibrium'', as well as the organizers of this very stimulating  trimester.

\end{document}